\documentclass[12pt]{amsart}
\usepackage{latexsym, amssymb, amsmath, amscd}
 \usepackage{eucal}

    \newtheorem{rema}{Remark}[section]
    \newtheorem{propo}[rema]{Proposition}

   \newtheorem{theo}[rema]{Theorem}
   \newtheorem{def-theo}[rema]{Definition-Theorem}
 \newtheorem{conj}[rema]{Conjecture}
   
    \newtheorem{lemma}[rema]{Lemma}
    \newtheorem{corol}[rema]{Corollary}
     
  \newtheorem{rmk}[rema]{Remark}

	\newcommand{\nno}{\nonumber}

 \newcommand{\pf}{{\it Proof:}\hspace{2ex}}
 
 \newcommand{\epfv}{\hspace{1em}$\Box$\vspace{1em}}

\newcommand{\cA}{{\mathcal A}}

\newcommand{\I}{{\operatorname{I}}}

\newcommand{\Ds}{{D_{\frak s}}}
\newcommand{\Dn}{{D_{\frak n}}}
\newcommand{\Ps}{{\phi_{\frak s}}}
\newcommand{\Pn}{{\phi_{\frak n}}}

\newcommand{\cD}{{\mathcal D}}

\newcommand{\cB}{{\mathcal B}}

\newcommand{\cE}{{\mathcal E}}

\newcommand{\Ker}{\operatorname{\rm Ker\,}}
\newcommand{\im}{\operatorname{Im}}

\newcommand{\Der}{{\cD er}}
\newcommand{\Eder}{{\cE der}}
\newcommand{\cEnd}{{\cE nd}}

\title[Idempotents in Intersection of the Kernel and the Image]
{Idempotents in Intersection of the Kernel and the Image of Locally Finite 
Derivations and $\cE$-derivations}

  \author{Wenhua Zhao}      
    \date{\today}
\address{Department of Mathematics, Illinois State University, Normal, IL 61761. Email: wzhao@ilstu.edu}

\begin{document}

\begin{abstract}
Let $K$ be a field of characteristic zero, $\cA$ a $K$-algebra and 
$\delta$ a $K$-derivation of $\cA$ or $K$-$\cE$-derivation of $\cA$ 
(i.e., $\delta=\operatorname{Id}_\cA-\phi$ for some $K$-algebra endomorphism $\phi$ 
of $\cA$). Motivated by the Idempotent conjecture proposed in \cite{Open-LFNED}, 
we first show that for all idempotent $e$ lying in both 
the kernel $\cA^\delta$ and the image $\im\delta\!:=\delta(\cA)$ of $\delta$, 
the principal ideal $(e)\subseteq \im \delta$ if $\delta$ is a locally finite $K$-derivation or a locally nilpotent $K$-$\cE$-derivation of $\cA$; and $e\cA,\, \cA e \subseteq \im\delta$ if $\delta$ is a locally finite $K$-$\cE$-derivation of $\cA$. 
Consequently, the Idempotent conjecture holds for all locally finite $K$-derivations and 
all locally nilpotent $K$-$\cE$-derivations of $\cA$. 
We then show that $1_\cA \in \im\delta$, (if and) only if $\delta$ is surjective,
which generalizes the same result \cite{GN, W} for locally nilpotent 
$K$-derivations of commutative $K$-algebras to locally finite 
$K$-derivations and $K$-$\cE$-derivations $\delta$ of all $K$-algebras $\cA$. 
%UDisk4/LFNED-Papers/Idems-in-Image-of-LFHD-Copy3.
\end{abstract}

   

%47B47 	Commutators, derivations, elementary operators, etc.

% 13G99	None of the above, but in this section of "Integral domains".

%16N40 Nil and nilpotent radicals, sets, ideals, rings
%
%16D70 Structure and classification (except as in 16Gxx), direct sum decomposition, cancellation

%16Wxx Rings and algebras with additional structure

%16W20 Automorphisms and endomorphisms

%16W25 Derivations, actions of Lie algebras

%16S34 Group rings [see also 20C05, 20C07], Laurent polynomial rings

%13N10: Rings of differential operators and their modules 
%14R15: Jacobian problem.
%33C45 Orthogonal polynomials and functions of hypergeometric type (Jacobi, Laguerre, Hermite, Askey scheme, etc.) [see also 42C05 for general orthogonal polynomials and functions]. 
%32C38 Sheaves of differential operators and their modules, $D$-modules. 
%32W99: Differential operators in several 
%    complex variables/None of the above, but in this section. 

\thanks{The author has been partially supported 
by the Simons Foundation grant 278638}

 \bibliographystyle{alpha}
    \maketitle

%\tableofcontents 

\renewcommand{\theequation}{\thesection.\arabic{equation}}
\renewcommand{\therema}{\thesection.\arabic{rema}}
\setcounter{equation}{0}
\setcounter{rema}{0}
\setcounter{section}{0}

\section{\bf Motivations and the Main Results} \label{S1}
Throughout the paper $K$ stands for a field of characteristic zero and $\cA$ for a $K$-algebra (not necessarily unital or commutative). We denote by $1_\cA$ or simply $1$ the identity element of $\cA$, if $\cA$ is unital, and 
$\I_\cA$ or simply $\I$ the identity map of $\cA$, if $\cA$ is clear in the context.  

A $K$-linear endomorphism $\eta$ of $\cA$ is said to be {\it locally nilpotent} (LN) 
if for each $a\in \cA$ there exists $m\ge 1$ such that $\eta^m(a)=0$, 
and {\it locally finite} (LF) if for each $a\in \cA$ the $K$-subspace spanned 
by $\eta^i(a)$ $(i\ge 0)$ is finite dimensional over $K$.   

By a $K$-derivation $D$ of $\cA$ we mean a $K$-linear map 
$D:\cA \to \cA$ that satisfies $D(ab)=D(a)b+aD(b)$ for all $a, b\in \cA$. 
By a $K$-$\cE$-derivation $\delta$ of $\cA$ 
we mean a $K$-linear map $\delta:\cA \to \cA$ such that for all 
$a, b\in \cA$ the following equation holds:
\begin{align}\label{ProdRule2}
\delta(ab)=\delta(a)b+a\delta(b)-\delta(a)\delta(b).  
\end{align}

It is easy to verify that $\delta$ is an $R$-$\cE$-derivation of $\cA$, 
if and only if $\delta=\I-\phi$ for some $R$-algebra endomorphism 
$\phi$ of $\cA$. Therefore an $R$-$\cE$-derivation is  
a special so-called $(s_1, s_2)$-derivation introduced 
by N. Jacobson \cite{J} and also a special 
semi-derivation introduced 
by J. Bergen in \cite{Bergen}. $R$-$\cE$-derivations have also been 
studied by many others under some different names such as 
$f$-derivations in \cite{E0, E} and 
$\phi$-derivations in \cite{BFF, BV}, etc..  

We denote by $\cEnd_K(\cA)$ the set of all 
$K$-algebra endomorphisms of $\cA$, $\Der_K(\cA)$ the set of all 
$K$-derivations of $\cA$, and $\Eder_K(\cA)$ the set of all 
$K$-$\cE$-derivations of $\cA$. Furthermore, for each $K$-linear endomorphism 
$\eta$ of $\cA$ we denote by $\im \eta$ the {\it image} of $\eta$, i.e., 
$\im\eta\!:=\eta(\cA)$, and $\Ker \eta$ the {\it kernel} of $\eta$.
When $\eta$ is an $R$-derivation or $R$-$\cE$-derivation, we also denote  
by $\cA^{{}^\eta}$ the kernel of $\eta$.  

It is conjectured in \cite{Open-LFNED} that the image  
of a LF $K$-derivation or $K$-$\cE$-derivation of $\cA$ 
possesses an algebraic structure, namely, 
a {\it Mathieu subspace}. The notion of Mathieu subspaces 
was introduced in \cite{GIC} and \cite{MS}, and is also called 
a {\it Mathieu-Zhao space} in the literature 
(e.g., see \cite{DEZ, EN, EH}, etc.) as first suggested 
by A. van den Essen \cite{E2}. 

The introduction of this new notion   
was mainly motivated by the study in \cite{Ma, IC} of the well-known 
Jacobian conjecture (see \cite{K, BCW, E}). 
See also \cite{DEZ}. But, a more interesting aspect 
of the notion is that it provides a natural  
generalization of the notion of ideals.

For some other studies on the algebraic structure of 
the image of a LF or LN $K$-derivations or 
$K$-$\cE$-derivations, see \cite{EWZ}, \cite{Open-LFNED}--\cite{OneVariableCase}.    

One motivation of this paper is the following so-called 
{\it Idempotent conjecture} proposed in \cite{Open-LFNED}, which is a weaker version 
of the conjecture mentioned above on the possible Mathieu subspace 
structure of the images of LF $K$-derivations and $K$-$\cE$-derivations.  

\begin{conj}\label{Idem-Conj}
Let $\delta$ be a LF (locally finite) $K$-derivation or a LF 
$K$-$\cE$-derivation of $\cA$ and $e\in \im \delta$ an idempotent of $\cA$, 
i.e., $e^2=e$. Then the principal (two-sided) ideal $(e)$ of $\cA$ 
generated by $e$ is contained in $\im\delta$.
\end{conj} 

Our first main result is the following theorem, which gives a partial 
positive answer to Conjecture \ref{Idem-Conj} above.

\begin{theo}\label{MainThm-1}
Let $K$ be a field of characteristic zero and $\cA$ a $K$-algebra 
(not necessarily unital). Then the following statements hold:
\begin{enumerate}
  \item[$1)$] for every locally finite $D\in \Der_K(\cA)$ and an idempotent 
$e\in \cA^D\cap\im D$, we have $(e)\subseteq\im D$;
\item[$2)$]  for every locally finite $\delta\in \Eder_K(\cA)$ and an idempotent 
$e\in \cA^\delta\cap\im\delta$, we have $e\cA,\,\cA e  \subseteq\im\delta$. Furthermore, 
if $\delta$ is locally nilpotent, we also have $(e) \subseteq\im\delta$. 
\end{enumerate}
\end{theo}

Note that for every $D\in \Der_K(\cA)$, it can be readily verify that 
all central idempotents  of $\cA$ lie in $\cA^D$. Furthermore, by Corollary \ref{LN-Central-Idem} that will be shown in Section \ref{S2},  this is also the case 
for every LN (locally nilpotent) $\delta\in \Eder_K(\cA)$. 
Therefore we immediately have the following 

\begin{corol}\label{Corol-1.1}
Assume that $\cA$ is a commutative $K$-algebra. Then 
Conjecture \ref{Idem-Conj} holds for all locally finite  
$D\in \Der_K(\cA)$ and all locally nilpotent $\delta\in \Eder_K(\cA)$.  
\end{corol}

For a different proof of the corollary above for commutative 
algebraic $K$-algebras, see \cite[Proposition 3.8]{Open-LFNED}. %\label{AlgCase2} 
%Note here that $\cA$ do NOT have to be algebraic. 
For a different proof of Theorem \ref{MainThm-1} for algebraic $K$-algebras (not necessarily commutative), 
see \cite[Corollary 3.9]{Open-LFNED}. 

Our second main result of this paper is the following 

\begin{propo}\label{Surj-one}
Assume that $\cA$ is unital and $\delta$ is a LF  
$K$-derivation or a LF $K$-$\cE$-derivation of $\cA$. 
Then $1_\cA \in \im\delta$, (if and) only if $\im\delta=\cA$, 
i.e., $\delta$ is surjective.   
\end{propo}

Two remarks about the proposition above are as follows. 

First, the proposition for LN $K$-derivations of commutative $K$-algebras was first proved 
by P. Gabriel and Y. Nouaz\'e \cite{GN} and later re-proved independently by D. Wright \cite{W}.
See also \cite{E}. 
During the preparation 
of this paper the author was informed that Arno van de Essen and Andrzej Nowicki 
have also proved the LF $K$-derivation case of the proposition 
for commutative $K$-algebras.  

Second, if $1_\cA \in \cA^\delta$, e.g., when $\delta\in \Der_K(\cA)$, 
the proposition follows immediately from Theorem \ref{MainThm-1}. But, 
if $1_\cA\not \in \cA^\delta$, the proof needs some other arguments  
(See Section \ref{S5}).   \\

%Our third main result is Theorem \ref{D-E-Corres} given in Section \ref{S2}, 
%which gives an one-to-one correspondence between the set 
%of all LN $K$-derivations of $\cA$ and the set of all LN 
%$K$-$\cE$-derivations of $\cA$. This theorem not only plays a crucial role in 
%the proofs of Theorem \ref{MainThm-1} and Proposition \ref{Surj-one}, 
%but also is interesting on its own right.  \\
 
{\bf Arrangement:} In Section \ref{S2} we recall and 
give a some shorter proof for van den Essen's one-to-one correspondence 
between the set of all LN $K$-derivations of $\cA$ and 
the set of all LN $K$-$\cE$-derivations of $\cA$ 
(See Theorem \ref{D-E-Corres}). We also derive some  
consequences of this important theorem that will be needed 
later in this paper. 
In Section \ref{S3} we show the $K$-derivation case 
of Theorem \ref{MainThm-1}. In Section \ref{S4} we show the $K$-$\cE$-derivation 
case of Theorem \ref{MainThm-1}. In Section \ref{S5} we give 
a proof for Proposition \ref{Surj-one}. \\

{\bf Acknowledgment:} The author is very grateful to Professors Arno van de Essen and Andrzej Nowicki for personal communications.

\section{\bf Van den Essen's One-to-One Correspondence between Locally Nilpotent Derivations 
and Locally Nilpotent $\cE$-Derivations} \label{S2}

Throughout this section, {\it $K$ stands for a field of characteristic zero, 
$\cA$ for a $K$-algebra (not necessarily unital or commutative) 
and $\I$ for the identity map of $\cA$.}

Denote by $\cD$ the set of all LN (locally nilpotent) $K$-derivations  
of $\cA$ and $\cE$ the set of all LN $K$-$\cE$-derivations of $\cA$. 
We define the following map:  
\begin{align}
\Xi: \cD&\to \cE \\
       D &\to \I-e^D, \nno
\end{align}
where $e^D\!:=\sum_{k=0}^\infty \frac{D^k}{k!}$.
 
With the setting as above we have the following remarkable one-to-one correspondence 
between $\cD$ and $\cE$, which was first proved by A. van den Essen in \cite{E0}. 
See also \cite[Proposition 2.1.3]{E}.   

\begin{theo}\label{D-E-Corres}
The map $\Xi:\cD\to\cE$ is an one-to-one correspondence between the sets $\cD$ and $\cE$ 
with the inverse map $\Xi^{-1}$ 
given by the following map:
\begin{align}
\Lambda: \cE&\to \cD \label{D-E-Corres-eq1} \\
       \delta &\to \ln(\I-\delta), \nno
\end{align}
where $\ln(\I-\delta)\!:=-\sum_{k=1}^\infty \frac{\delta^k}k$.   

%$2)$ For all $D\in \cD$, let $\delta=\Xi (D)$. Then $\cA^D=\cA^\delta$ and $\im D=\im \delta$.  
\end{theo}

For the sake of completeness, we here give a proof for the theorem above, 
which is some shorter than the one given in \cite[Proposition 2.1.3]{E}. 

%\begin{lemma}
%$1)$ Let $E$ be a locally nilpotent $K$-linear endomorphism of $\cA$. 
%Then $\I-E$ is invertible.
%
%$2)$ Let $f(t)$ be a nonzero formal power series in $t$ and $U$, $V$ LN $K$-linear endomorphism of $\cA$. If $f(U)=f(V)$, then $U=V$. 
%\end{lemma}
%
%\pf $1)$ Let $F=\sum_{i=0}^\infty E^i$. Since $E$ is LN, $F$ is well-defined $K$-linear map of $\cA$. It is easy to see that $F=E^{-1}$.
%
%$2)$ Let $u\in \cA$. Since $U$ and $V$ are LN, there exists $n\ge 1$ 
%such that $U^n(u)=0=V^n(u)$. We use the (backwards) induction to 
%show that $U^i(u)=V^i(u)$ for all $1\le i\le n$.
%
%Assume $U^j(u)=V^j(u)$ for some $2\le j\le n$. Write $f(
%\epfv
%
%\begin{lemma}
%
%\end{lemma}

First, the following lemma can be easily verified by induction, 
as noticed in \cite{E0, E}. 

\begin{lemma}\label{n-E-Lebniz}
 Let $\cB$ be a ring and $\delta$ an $\cE$-derivation of $\cB$. Then for all 
$a, b\in \cB$ and $n\ge 1$, we have
\begin{align}\label{n-E-Lebniz-eq1}
\delta^n(ab)=\sum_{i=0}^n \binom ni \delta^i(a)\, \delta^{n-i}(\I-\delta)^i(b).
\end{align}
 \end{lemma}

%\pf We use induction on $n$. When $n=1$, Eq.\,(\ref{n-E-Lebniz-eq1}) becomes 
%\begin{align}\label{n-E-Lebniz-peq1}
%\delta (ab) =a\delta (b)+\delta(a)(\I-\delta)(b),
%\end{align}
%which is the same as Eq.\,(\ref{ProdRule2}). 
%
%So we assume that Eq.\,(\ref{n-E-Lebniz-eq1}) holds for some $n\ge 1$. Consider
%\begin{align*}
%\delta^{n+1}(ab) =\delta \big(\delta^n(ab)\big)= \sum_{i=0}^n \binom ni \Lambda_{n i}(a, b),
%\end{align*}
%where by the induction assumption and 
%Eq.\,(\ref{n-E-Lebniz-peq1}), $\Lambda_{n i}(a, b)$ is equal to 
%\begin{align*} 
%% & \delta^i(a)  \delta^{n-i+1}(\I-\delta)^i(b)+ \left(
%%\delta^{i+1}(a) \delta^{n-i}(\I-\delta)^i(b)-\delta^{i+1}(a)  \delta^{n-i+1}(\I-\delta)^i(b)\right) \\
%&\delta^i(a) \delta^{n-i+1}(\I-\delta)^i(b)+ \delta^{i+1}(a) \left( 
%(\I-\delta)\delta^{n-i}(\I-\delta)^i\right) (b) \\
%&=\delta^i(a) \delta^{n-i+1}(\I-\delta)^i(b)+\delta^{i+1}(a) \delta^{n-i}(\I-\delta)^{i+1}(b).
%\end{align*}
%Therefore we see that $\delta^{n+1}(ab)$ is equal to 
%\begin{align*}
%&\sum_{i=0}^n \binom ni  \delta^i(a) \delta^{n-i+1}(\I-\delta)^i(b)+
% \sum_{i=0}^n \binom ni \delta^{i+1}(a) \delta^{n-i}(\I-\delta)^{i+1}(b)\\
%&= \sum_{i=0}^{n+1} \left(\binom {n}i+\binom n{i-1}\right)  \delta^i(a) \delta^{n-i+1}(\I-\delta)^i(b), \\
%\intertext{where $\binom n{-1}=\binom n{n+1}=0$.}  
%&= \sum_{i=0}^{n+1} \binom {n+1}i  \delta^i(a) \delta^{n-i+1}(\I-\delta)^i(b).
%\end{align*}
%Hence by induction the lemma follows.
%\epfv

Now we can show Theorem \ref{D-E-Corres} as follows.\\

\underline{\bf Proof of Theorem \ref{D-E-Corres}:} \, 
First, since $D$ is LN, $e^D$ is well-defined.  
It is well-known (and also easy to check directly) 
that $e^D$ is a $K$-algebra automorphism of $\cA$. 
Hence $\Xi(D)$ is a $K$-$\cE$-derivation of $\cA$.
Consider
\begin{align}\label{D-E-Corres-peq1}
\Xi(D)=\I-e^D= -\sum_{n=1}^\infty \frac{D^n}{n!}=D h(D)=h(D)D,  
\end{align} 
where 
\begin{align}\label{D-E-Corres-peq1.2}
h(D) =-\I-\sum_{n=2}^\infty \frac{D^{n-1}}{n!}=-\I-
D\sum_{n=2}^\infty \frac{D^{n-2}}{n!}.  
\end{align}

Since $D$ is LN, and $D$ and $h(D)$ commute, by Eq.\,(\ref{D-E-Corres-peq1}) 
$\Xi(D)$ is also LN. Therefore, $\Xi$ is indeed a map from $\cD$ to $\cE$.

%Furthermore, by Lemma \ref{InvLma} with 
%$F=-\I$ and $G=D\sum_{n=2}^\infty \frac{D^{n-2}}{n!}$ 
%we see that $h(D)$ is an $K$-linear automorphism of 
%$\cA$ that commutate with both $D$ and $\Xi(D)$. 
%Hence statement $2)$ follows also from 
%Eq.\,(\ref{D-E-Corres-peq1}).

Next, we show that $\Lambda (\delta)\in \cD$ for all 
$\delta\in \cE$.  
Set $D_\delta\!:=\Lambda(\delta)$. Then  
\begin{align} \label{D-E-Corres-peq2}
D_\delta\!=\ln(\I-\delta)=
-\sum_{n=1}^\infty \frac{\delta^n}{n}=\delta g(\delta)=g(\delta)\delta,
\end{align}
where 
\begin{align} \label{D-E-Corres-peq2.2}
g(\delta)=-\sum_{n=1}^\infty \frac{\delta^{n-1}}{n}=
-\I-\delta \sum_{n=2}^\infty \frac{\delta^{n-2}}{n}.
\end{align}

Since $\delta$ is LN, and $\delta$ and $g(\delta)$ commute, 
by Eq.\,(\ref{D-E-Corres-peq2}) $D_\delta$ is also LN.  
 
Now, let $x, y\in \cA$. Then by Lemma \ref{n-E-Lebniz} we have  
\begin{align}
D_\delta(xy)&=-\sum_{n=1}^\infty \frac1n \delta^n(xy) \label{D-E-Corres-peq3}\\
%\intertext{Applying Lemma \ref{n-E-Lebniz}}  
 &=-\sum_{n=1}^\infty \frac1n \sum_{i=0}^n\binom{n}i \delta^i(x) \delta^{n-i}(1-\delta)^i(y)\nno \\
&=-\sum_{i=0}^\infty \delta^i(x) S_i(y), \nno 
\end{align}
where for each $i\ge 0$, 
\begin{align}\label{D-E-Corres-peq4}
 S_i(y) =\sum_{n=i}^\infty \frac1n \binom{n}i \delta^{n-i}(1-\delta)^i(y) 
=(1-\delta)^i \sum_{n=i}^\infty \frac1n \binom{n}i \delta^{n-i}(y). 
\end{align}
In particular, by Eq.\,(\ref{D-E-Corres-peq2}) and the equation above we have 
\begin{align}\label{D-E-Corres-peq4b}
S_0(y)=-D_\delta(y). 
\end{align}

\vspace{3mm}

{\bf Claim}: {\it $S_i(y)=\frac1i y$ for all $i\ge 1$.} \\

\underline{\it Proof of Claim}: For each $i\ge 1$ 
we introduce the formal power series  
\begin{align}
f_i(t)\!:&=(1-t)^i \sum_{n=i}^\infty \frac1n \binom{n}i t^{n-i} \label{D-E-Corres-peq5} \\
&=\frac{(1-t)^i}{i!} \sum_{n=i}^\infty (n-1)(n-2)\cdots (n-i+1)t^{n-i}. \nno
\end{align}
Then $S_i(y)=f_i(\delta)(y)$. On the other hand, 
we have the following identity of formal power series:
$$ 
 \frac{(i-1)!}{(1-t)^i}=-\frac{d^i}{dt^i} \ln (1-t)= \sum_{n=i}^\infty (n-1)(n-2)\cdots (n-i+1)t^{n-i}
$$

By Eq.\,(\ref{D-E-Corres-peq4}) and the identity above  
  we have $f_i(t) =1/i$. Hence  
$S_i(y)=f_i(\delta)(y)=\frac{1}{i} y$ and 
the claim follows. \\

Now by Eqs.\,(\ref{D-E-Corres-peq2}), (\ref{D-E-Corres-peq3}), (\ref{D-E-Corres-peq4b}) and the claim 
above we have 
\begin{align*}
D_\delta (xy)&=-\sum_{i=0}^\infty \delta^i(x) S_i(y) 
=-xS_0(y)- \sum_{i=1}^\infty  \delta^i(x) S_i(y) \\&=
xD_\delta(y)- \sum_{i=1}^\infty \frac1i \delta^i(x)y
=xD_\delta(y)+D_\delta(x)y.
\end{align*}
Therefore $\Lambda(\delta)=D_\delta$ is a LN   
$K$-derivation of $\cA$, i.e., $\Lambda$ is indeed 
a map from $\cE$ to $\cD$.  Since $\Xi$ and 
$\Lambda$ are obviously inverse to each other, 
we see that $\Xi$  gives an one-to-one 
correspondence between $\cD$ to $\cE$, 
i.e., the theorem follows. 
\epfv

Next, we derive some consequences of Theorem \ref{D-E-Corres}.  
But, we first need to show the following lemma. Although it is 
almost trivial, it will be frequently used throughout 
the rest of the paper.  

\begin{lemma}\label{InvLma}
Let $R$ be a ring and $\cB$ an $R$-algebra.
Let $F$ and $G$ be two commuting $R$-linear endomorphisms 
of $\cB$ such that $F$ is invertible and $G$ is LN (locally nilpotent). 
Then $F-G$ is an $R$-linear automorphism 
with the inverse map given by 
$$
(F-G)^{-1}=\sum_{k=0}^\infty G^kF^{-k-1}.
$$
\end{lemma}

\pf Note that $F-G=(\I-GF^{-1})F$.
Since $F$ commutes with $G$, so does $F^{-1}$. 
Hence $U\!:=GF^{-1}$ is LN, 
for $G$ is LN. Therefore the formal power series 
$\sum_{k=0}^\infty U^k$ is a well-defined $R$-linear 
endomorphism of $\cA$, which gives 
the inverse map of $\I-U$. Hence, $F^{-1}\sum_{k=0}^\infty U^k$ gives 
the inverse map of $F-G$, from which the 
lemma follows.
\epfv

\begin{corol}\label{D-E-Corres-Corol}
Let $\cD$, $\Xi$ be as in Theorem \ref{D-E-Corres}, and $D\in\cD$. 
Set $\delta=\Xi(D)$. Then $\cA^D=\cA^\delta$ and $\im D=\im \delta$. 
\end{corol} 

\pf First, since $\delta=\Xi(D)$, we have $D=\Lambda(\delta)=:\!D_\delta$ by Theorem \ref{D-E-Corres}. 
Second, by Lemma \ref{InvLma} with $F=-\I$ and 
$G=D\sum_{n=2}^\infty \frac{D^{n-2}}{n!}$,   
the $K$-linear map $h(D)$ in Eq.\,(\ref{D-E-Corres-peq1.2}) 
is a $K$-linear automorphism of $\cA$. 
Therefore, we have $\cA^D=\cA^\delta$ by Eq.\,(\ref{D-E-Corres-peq1}). 
Furthermore, we also have $\im \delta \subseteq \im D$ by Eq.\,(\ref{D-E-Corres-peq1}),  
and $\im D \subseteq \im \delta$ 
by Eq.\,(\ref{D-E-Corres-peq2}), 
whence $\im D = \im \delta$, and the corollary 
follows.
\epfv

\begin{corol} \label{LN-Central-Idem}
Let $\delta$  be an arbitrary $K$-derivation of $\cA$ or a LN   
$K$-$\cE$-derivation of $\cA$. Then all central idempotents of $\cA$ 
lie in $\cA^\delta$.
\end{corol} 

\pf If $\delta$ is a $K$-derivation of $D$, the corollary actually 
holds regardless of the characteristic of $K$, which 
can be seen as follows. 
Since $De=De^2=2eDe$, we have $(1-2e)De=0$. 
Since $(1-2e)^2=1-4e+4e^2=1$, $1-2e$ is a unit of $\cA$. Hence $De=0$.

If $\delta$ is a LN $K$-$\cE$-derivation of $\cA$, then by Theorem \ref{D-E-Corres}, 
$\delta=\I-e^D$ for some LN $K$-derivation $D$ of 
$\cA$. Since $De=0$ as shown above, we have $\delta e=0$. 
\epfv

%Furthermore, it is easy to see from Theorem \ref{D-E-Corres}   
%and its proof that the following  corollary  
%also holds. 
% 
%%\begin{corol}
%%A $K$-algebra $\cA$ has no nonzero LN $K$-derivations, if and only if 
%%$\cA$ has no nonzero LN $K$-$\cE$-derivations.
%%\end{corol}
%
%\begin{corol}
%$1)$ For all $D\in \cD$ we have $\Xi(-D)=\I-(\I-\delta)^{-1}$, where 
%$\delta=\Xi(D)$. 
% 
%$2)$ For all $u\in \cA$, $D\in \cD$ and $k\ge 1$, let $\delta=\Xi(D)$. Then  
%$D^k(u)=0$, if and only if $\delta^k(u)=0$.
% 
%$3)$  $\Xi$ maps nilpotent $K$-derivations of $\cA$ to nilpotent 
%$K$-$\cE$-derivations of $\cA$, and similarly for the corresponding $\Lambda=\Xi^{-1}$.
%
%$4)$  $\Xi$ maps commuting LN $K$-derivations of $\cA$ to commuting LN $K$-$\cE$-derivations  
%of $\cA$, and similarly for the corresponding $\Lambda=\Xi^{-1}$. 
%\end{corol} 

\section{\bf The Derivation Case of Theorem \ref{MainThm-1}} \label{S3}

In this section we give a proof of Theorem \ref{MainThm-1} for 
LF (locally finite) $K$-derivations. 
%
%and $K$-$\cE$-derivations. 
%By Theorem \ref{D-E-Corres} it suffices to show the case of 
%locally nilpotent $K$-derivations. 
%
Throughout this section we let $K$ and $\cA$ be as in  
Theorem \ref{MainThm-1}, $D$ a LF $K$-derivation of $\cA$,  
and $e$ an idempotent in $\cA^D\cap\im D$. 

Let $s\in \cA$ such that $Ds=e$. Since $De=0$, 
we have $D(ese)=e(Ds)e=e$. So replacing $s$ by 
$ese$ we assume $s\in e\cA e$. 
Furthermore, for convenience we set $\mathbf{s^0=e}$. 
Then with the setting above it is easy 
to see that for all $i, k\ge 0$, we have 
$es^i=s^ie=s^i$ and 
\begin{align}\label{SpecialPower}
D^i(s^k)=\begin{cases}
k(k-1)\cdots(k-i+1)s^{k-i} &\text{ if } i\le k \\
0 &\text{ if } i>k.
\end{cases}
\end{align}
 
We first consider the case that $D$ is LN (locally nilpotent).

\begin{lemma}\label{phi-psi-Lma}
Assume that $D$ is LN. For all $a\in \cA$ set
\begin{align}
\phi_{-s}(a)&=\sum_{i=0}^\infty \frac{(-1)^i}{i!} D^i(a)s^i \label{Def-phi}\\
\psi_{-s}(a)&=\sum_{i=0}^\infty \frac{(-1)^i}{i!} s^i D^i(a).
\end{align}
Then  
$\phi_{-s}(a), \psi_{-s}(a)\in \cA^D$ and 
\begin{align}
ae=\sum_{j=0}^\infty \frac1{j!} \phi_{-s}\big(D^j(a)\big) \, s^j \label{phi-psi-Lma-eq1}\\
ea=\sum_{j=0}^\infty \frac1{j!} s^j \psi_{-s}\big(D^j(a)\big). \label{phi-psi-Lma-eq2}
\end{align}
\end{lemma}

Note that the case when $\cA$ is commutative and $e=1$ the lemma has been proven 
in \cite{GN, W}. See also \cite{E}. The main idea of the proof given below is to modify the proof in  \cite{GN, W, E} to the more general case in the lemma. \\
%Then by a similar idea as in the proof of Corollary 1.3.23 in \cite{E}, 
%which is the case when $\cA$ is commutative and $e=1$, 
%we can prove the following more general lemma.  

\pf First, by Eqs.\,(\ref{SpecialPower}) and (\ref{Def-phi}) we have 
\begin{align*}
D\phi_{-s}(a)&= \sum_{i=0}^\infty \frac{(-1)^i}{i!} D
\big( D^i(a)s^i \big ) \\&
 =\sum_{i=0}^\infty \frac{(-1)^i}{i!} 
\big( D^{i+1}(a) s^i  + iD^{i}(a) s^{i-1}\big ) \\
&=\sum_{j=0}^\infty \frac{(-1)^j+(-1)^{j+1}}{j!} D^{j+1}(a)s^j   
 =0.
\end{align*}

Hence $\phi_{-s}(a)\in \cA^D$. 
The proof of $\psi_{-s}(a)\in \cA^D$ 
is similar. 

Next we show Eq.\,(\ref{phi-psi-Lma-eq1}). 
The proof of Eq.\,(\ref{phi-psi-Lma-eq2}) is similar.
\begin{align*}
\sum_{j=0}^\infty \frac1{j!} \phi_{-s}\big(D^j(a)\big) \, s^j&=
\sum_{j=0}^\infty \frac1{j!} \sum_{i=0}^\infty \frac{(-1)^i}{i!}  
D^{i+j}(a)s^{i+j} \\ 
&=ae+\sum_{n=1}^\infty n! \left( \sum_{\substack{i+j=n\\i,\,j\ge 0}} 
\binom ni (-1)^i\right) D^n(a) s^n\\
&=ae.  
\end{align*}
\epfv

\begin{propo}
Let $\delta$ be a LN (locally nilpotent) $K$-derivation or $K$-$\cE$-derivation 
of $\cA$ and $e\in \cA^\delta\cap \im \delta$ a nonzero idempotent. 
Let $s\in \cA$ such that $\delta s=e$. Replacing $s$ by $ese$ 
we assume $s\in e\cA e$. Set ${\bf s^0=e}$. Then we have
\begin{enumerate}
\item[$1)$] if $\sum_{i=0}^n c_i s^i=0$ or $\sum_{i=0}^n s^ic_i=0$ with $c_i\in \cA^\delta$, 
then $c_ie=0$ for all $0\le i\le n$. In particular, 
$s$ is transcendental over the field $Ke$;
\item[$2)$] $\cA e=\cA^\delta[s]$ and $e\cA =[s]\cA^\delta$, 
where $\cA^\delta[s]$ (resp., $[s]\cA^\delta$) is the $K$-algebra of 
all polynomials $f(s)$ of the form 
$f(s)=\sum_{i\ge 0}a_is^i$ (resp., $f(s)=\sum_{i\ge 0}s^ia_i$) 
with $a_i's$ in $\cA^\delta$; 
\item[$3)$] if $\delta\in \Der_K(\cA)$, then $\delta\,\mid_{e\cA}=\frac{d}{ds}$ and 
$\delta\,\mid_{\cA e}=\frac{d}{ds}$;
\item[$4)$] if $\delta\in \Eder_K(\cA)$, then $\delta\,\mid_{e\cA}=\I_{e\cA}-\phi$ 
  (resp., $\delta\,\mid_{\cA e}=\I_{\cA e}-\psi$), where $\phi$ (resp., $\psi$) is the $\cA^\delta$-algebra endomorphism of $e\cA$ (resp., $\cA e$) that maps $s$ to $s+e$.   
\end{enumerate}
\end{propo}

The proof of the $K$-derivation case of the proposition above is similar as 
the proof of \cite[Proposition $1.3.21$]{E}. The $K$-$\cE$-derivation case follows 
from Theorem \ref{D-E-Corres} and the $K$-derivation case of the proposition.    
So we skip the detailed proof of this proposition here.

% 
%Note that the theorem is not true if char.\,$K=p>0$.
%
%\begin{exam}
%Let $\cA=\bF_p[x]/(x^p)$ and $D=\p_x$. Note that $D$ preserves 
%the ideal $(x^p)$, and hence induces a derivation of $\cA$ 
%(still denote it by $D$).
%Then $D(x)=1$ but $D$ is {\bf NOT} surjective 
%(for $x^{p-1}\not \in \im D$). 
%
%This gives also an example for an algebraic 
%derivation whose image contains an non-zero idempotent. 
%\end{exam}

%\begin{corol}
%Let $\cA$ be an algebraic $\bQ$-algebra and  
%$D$ a $\bQ$-LND. Then $1\not\in\im D$. 
%\end{corol}
%
%\pf Assume that $Ds=1$. Then by Theorem \ref{}, 
%$s$ is transcendental over $K$ for $K\subseteq \cA^D$. 
%Contradiction.
%\epfv
%
%Since each finite dimensional $K$-algebra is algebraic over $K$, 
%we immediately have 
%
%\begin{corol}
%Let $K$ be a field of characteristic zero and 
%$\cA$ be a finite dimensional $K$-algebra. 
%Then for each $K$-LND of $\cA$ we have $1\not \in \im D$, 
%i.e., $D$ has no slice.
%\end{corol}
% 

%Now we need to see if $\im D$ can contains some 
%idempotents.
%
%\begin{propo}
%Let $\cA$ be a $\bQ$-algebra, $D$ a $\bQ$-LND and $e$ 
%an central idempotent of $\cA$. Assume that $e\in \im D$, 
%then the principal ideal $(e)=e\cA\subseteq \im D$. 
%\end{propo}

From the proposition above we also have the following 

\begin{corol}
Assume further that $\cA$ is algebraic over $K$.  
Let $\delta$ be a LN $K$-derivation or $K$-$\cE$-derivation 
of $\cA$. Then $\cA^\delta\cap \im \delta$ does not contain any nonzero 
idempotent of $\cA$. 
\end{corol}

Actually, the LN condition on $\delta$ in the corollary above can be dropped.
See \cite[Corollary 3.9]{Open-LFNED}. For more results on the idempotents in the image   
of LF or LN $K$-derivations and $K$-$\cE$-derivations 
of algebraic $K$-algebras, see \cite{Alg-LFNED} and \cite{Open-LFNED}. 

Next, we consider Theorem \ref{MainThm-1}, first,  
for all LN $K$-derivations and  
$K$-$\cE$-derivations of $\cA$.

\begin{lemma}\label{LN-MainThm-1}
Theorem \ref{MainThm-1} holds for all LN $K$-derivations and  
$K$-$\cE$-derivations of $\cA$.
\end{lemma}

\pf Note first that by Theorem \ref{D-E-Corres} and Corollary \ref{D-E-Corres-Corol},   
it suffices to show the lemma for all LN $D\in \Der_K(\cA)$. Let $e$, $s$ and $a$ 
be as in Lemma \ref{phi-psi-Lma}. Then $\phi_{-s}(a), \psi_{-s}(a)\in \cA^D$ 
for all $a\in \cA$. By Eq.\,(\ref{phi-psi-Lma-eq1}) we see that $D$ maps $\sum_{j=0}^\infty \frac1{(j+1)!} \phi_{-s}\big(D^j(a)\big) \, s^{j+1}$
to $ae$, and by Eq.\,(\ref{phi-psi-Lma-eq2})   $D$ maps 
$\sum_{j=0}^\infty \frac1{(j+1)!} s^{j+1} \psi_{-s}\big(D^j(a)\big)$
to $ea$. Hence $ea, ae\in \im D$ for all $a\in \cA$.

To show $aeb\in \im D$ for all $a, b\in \cA$, note first that 
by Eq.\,(\ref{phi-psi-Lma-eq1}) for $ae$ and Eq.\,(\ref{phi-psi-Lma-eq2}) for $eb$ 
we have 
$$
aeb=\sum_{i, j=0}^\infty \frac1{i!j!} \phi_{-s}\big(D^i(a)\big) \, s^{i+j} \psi_{-s}\big(D^j(b)\big). 
$$
Then $D$ maps $\sum_{i, j=0}^\infty \frac1{i!j!} \phi_{-s}\big(D^i(a)\big) \, 
\frac{s^{i+j+1}}{i+j+1} \psi_{-s}\big(D^j(b)\big)$ to $aeb$. Therefore, 
we have  $(e)\subseteq D$, i.e., 
Theorem \ref{MainThm-1}, $1)$ holds for $D$, 
as desired.
\epfv

Now we assume that $D$ is LF and consider the case 
that the base field $K$ is algebraically closed. 
In this case $D$ has the Jordan-Chevalley decomposition 
$D=\Dn+\Ds$ over $K$ (see \cite[Proposition $1.3.8$]{E}) 
%or \cite[Proposition $4.2$]{H}) 
such that $\Ds$ is semi-simple and $\Dn$ is LN. 
%Furthermore, by \cite[Proposition $1.3.13$], \cite{E1}) 
%both $\Ds$ and $\Dn$ are $K$-derivations.  

Let $\Lambda$ be the set of all distinct 
eigenvalues of $\Ds$ and $\cA_\lambda$ 
$(\lambda\in \Lambda)$ the corresponding 
eigenspace $\Ds$.  
Then $\cA$ has the following direct 
sum decomposition:
\begin{align}
\cA=\oplus_{\lambda\in \Lambda} \cA_\lambda.\label{D-Decomp}
\end{align}

Actually, the decomposition above gives a $K$-algebra grading of $\cA$, i.e., $\cA_\lambda\cA_\mu\subseteq \cA_{\lambda+\mu}$ for all 
$\lambda, \mu\in \Lambda$. This is because 
$\Ds$ and $\Dn$  by \cite[Proposition 1.3.13]{E} are also   
$K$-derivations of $\cA$. % (e.g., see \cite[Proposition 1.3.13]{E}).  
In particular, $\cA_0$ is a
$K$-subalgebra of $\cA$. Furthermore, each $\cA_\lambda$ is $D$ 
(and also $\Ds$ and $\Dn$) invariant. Therefore we have  
\begin{align}
\im D=\oplus_{\lambda\in \Lambda} D(\cA_\lambda). \label{ImD-Decomp}
\end{align}

\begin{lemma} \label{D-preLma(-1)}
Assume that $K$ is algebraically closed. Then
\begin{enumerate}
\item[$1)$] $\cA^D\subseteq \cA_0$.
\item[$2)$] $\im D=\Dn(\cA_0) \oplus 
\bigoplus_{0\ne \lambda\in \Lambda} \cA_\lambda$. 
\item[$3)$] $(e) \subseteq \im D$ for all idempotents  
 $e\in \cA^D\cap\im D$.  
\end{enumerate}
\end{lemma}

\pf $1)$ By \cite[Proposition $1.3.9$, $i)$]{E} we have 
\begin{align} \label{D-preLma(-1)-peq1}
\cA^D=\Ker \Ds\cap \Ker \Dn.
\end{align}
Since $\Ker \Ds=\cA_0$, we hence 
have $\cA^D\subseteq \cA_0$. 

$2)$ For each $0\ne \lambda\in \Lambda$, we have 
$D(\cA_\lambda)\subseteq \cA_\lambda$ 
and 
$$
D \mid_{\cA_\lambda}=\Ds\mid_{\cA_\lambda}+\Dn\mid_{\cA_\lambda}
=\lambda \I_{\cA_\lambda}+\Dn\mid_{\cA_\lambda}.
$$ 
Since $\Dn$ is LN,   
$D\mid_{\cA_\lambda}$ by Lemma \ref{InvLma} is a $K$-linear automorphism 
of $\cA_\lambda$, whence $\cA_\lambda\subseteq \im D$ for all 
$0\ne \lambda\in \Lambda$. Note that $D\mid_{\cA_0}=\Dn\mid_{\cA_0}$. 
Then by Eq.\,(\ref{ImD-Decomp}) the statement follows.  

$3)$ Note that $\Dn$ is a LN $K$-derivation of $\cA$ 
(as pointed out above)   
and $e\in \Ker \Dn$ by Eq.\,(\ref{D-preLma(-1)-peq1}). 
Applying lemma \ref{LN-MainThm-1} to $\Dn$ we have 
$(e) \subseteq \im \Dn$. 
Therefore it suffices to show $\im\Dn\subseteq \im D$. 

Since $\cA_\lambda$ is $\Dn$ invariant for all $\lambda\in \Lambda$, 
we have 
$$
\im \Dn=\Dn(\cA_0)\oplus \bigoplus_{0\ne \lambda\in\Lambda} \Dn (\cA_\lambda).
$$ 
Since $\Dn (\cA_\lambda)\subseteq \cA_\lambda$ for all $\lambda\in \Lambda$, 
by statement $2)$ we have $\im\Dn\subseteq \im D$, as desired. 
 \epfv
%
%$3)$ By statement $1)$ $e\in \cA_0$, and then   
%by Statement $2)$ we have $e\in \Dn(\cA_0)$. 
%Since $\Dn\mid_{\cA_0}=D\mid_{\cA_0}$ is a LN $K$-derivation 
%of $\cA_0$, applying lemma \ref{LN-MainThm-1} to the $K$-algebra $\cA_0$ and 
%$\Dn\,\mid_{\cA_0}$ we have $e\cA_0, \, \cA_0 e,\, \cA_0 e\cA_0 
%\subseteq \Dn(\cA_0)=D(\cA_0) \subseteq \im D$.  
%Since $e\in \cA_0$, we also have 
%$e\cA_\lambda, \, \cA_\lambda e \subseteq \cA_\lambda$ 
%and $\cA_\lambda e\cA_\mu \subseteq \cA_{\lambda+\mu}$ 
%for all $\lambda, \mu \in \Lambda$. Then by statement $2)$ it suffices to 
%show $\cA_\lambda e\cA_{-\lambda} \subseteq \im D$ 
%for all $0\ne \lambda\in \Lambda$. 
%
%Let $0\ne \lambda\in \Lambda$ and $s\in \cA_0$ such that $Ds=e$. Then for all $u_\lambda\in \cA_\lambda$ and $u_{-\lambda}\in \cA_{-\lambda}$, we have 
%\begin{align}
%\Ds(u_\lambda s u_{-\lambda})& = (\Ds u_\lambda) s u_{-\lambda}+u_\lambda (\Ds s) u_{-\lambda}+ u_\lambda e (\Ds u_{-\lambda}) \\
%&=\lambda u_\lambda s u_{-\lambda}+u_\lambda e u_{-\lambda}-\lambda u_\lambda s u_{-\lambda}\nno \\
%&=u_\lambda e u_{-\lambda}. \nno
%\end{align}
%
%we see that statement $3)$ follows.
%\epfv

Now we drop the assumption that $K$ is algebraically closed and show 
Theorem \ref{MainThm-1} for all LF $K$-derivations $D$ of $\cA$. \\

\underline{\bf Proof of Theorem \ref{MainThm-1}, 1)}: 
Let $\bar K$ be the algebraic closure of $K$ and $\bar\cA=\bar K\otimes_K\cA$.
Then we may identify $\cA$ as a $K$-subalgebra of $\bar \cA$ in the standard way.  
Denote by $\bar D$ the $\bar K$-linear extension map of $D$ from $\bar \cA$ to $\bar \cA$. 
Then it is easy to see that $\bar D$ is a LF  
$\bar K$-derivation of $\bar\cA$.  

Let $e$ be an idempotent in $\cA^D\cap \im D$. Then $e$ also lies 
in $ \bar\cA^{\bar D}\cap \im \bar D$, and 
by Lemma \ref{D-preLma(-1)}, $3)$ we have   
$e\bar \cA, \, \bar\cA e,\, \bar\cA e\bar\cA \subseteq \bar D(\bar\cA)$,    
whence $e\cA,\, \cA e,\, \cA e\cA \subseteq \cA \cap \bar D(\bar\cA)$. 
Then Theorem \ref{MainThm-1}, $1)$ immediately follows  
from the lemma below.
\epfv

\begin{lemma}\label{lin-bar-Lma}
Let $V$ be a vector space over a field $K$ (not necessarily of characteristic zero) 
and $f$ a $K$-linear endomorphism of $V$. Let $L$ be a field extension of $K$, 
$\bar V\!:=L\otimes_K V$ and $\bar f$ the $L$-linear extension map of $f$ from 
$\bar V$ to $\bar V$. Identify $V$ as a $K$-subspace of $\bar V$ 
in the standard way. Then $f(V)=V \cap \bar f(\bar V)$. 
\end{lemma}

\pf It suffices to show that for each $v\in V$, 
there exists $u\in \cA$ such that $f(u)=v$  if (and only if) there exists 
$\bar u\in \bar \cA$ such that $\bar f(\bar u)=v$. 
Let $v_i$ $(1\le i\le n)$ be the $K$-linearly independent vectors in $V$ 
such that $v\in \text{Span}_K\{v_i\,|\, 1\le i\le n\}$ 
and $\bar u \in \text{Span}_L \{v_i\,|\, 1\le i\le n\}$. 
By using the coordinates of $\bar u$ and $v$, and the transformation  
matrix of $f$ with respect to $\{v_i\,|\, 1\le i\le n\}$, we see that 
the problem becomes the following problem on linear systems: 
{\it for all $y\in K^n$ and 
$n\times n$ matrix $A$ with entries in $K$, the linear system $Ax=y$ 
has a solution $x$ in $K^n$ if (and only if) it has a solution in 
$L^n$}. But this can be easily verified, 
e.g., by applying elementary row operations 
to transform $A$ into an up-triangular matrix.  
\epfv

\section{\bf The $\cE$-Derivation Case of Theorem \ref{MainThm-1}}\label{S4}

Throughout this section we let $K$ and $\cA$ be as in Theorem \ref{MainThm-1} and 
fix a LF (locally finite) $K$-$\cE$-derivation $\delta$ of $\cA$. 
Write $\delta=\I-\phi$ for some $K$-algebra endomorphism 
$\phi$ of $\cA$. Note that $\cA^\delta=\cA^\phi\!:=\{ u\in \cA\,|\, \phi(u)=u\}$ 
and $\phi$ is also LF.  

%Furthermore, by Theorem \ref{D-E-Corres} 
%and the $K$-derivation case of Theorem \ref{MainThm-1} proved in the previous section  
%we immediately have the following
%
%\begin{lemma}\label{Lma-4.1}
%If $\delta$ is LN, then 
%Theorem \ref{MainThm-1} holds for $\delta$. 
%\end{lemma}

We first assume that $K$ is algebraically closed. 
In this case 
$\phi$ has the Jordan-Chevalley decomposition $\phi=\Pn+\Ps$ 
over $K$ (e.g., see \cite[Proposition $1.3.8$]{E}) 
%or \cite[Proposition $4.2$]{H}) 
such that $\Ps$ 
is semi-simple and $\Pn$ is LN (locally nilpotent). 
%Furthermore, by \cite[Proposition $1.3.13$], \cite{E1}) 
%both $\Ds$ and $\Dn$ are $K$-derivations.  

Let $\Lambda$ be the set of all distinct 
eigenvalues of $\Ps$ and $\cA_\lambda$ 
$(\lambda\in \Lambda)$ the corresponding 
eigenspace of $\Ps$.  
Then $\cA$ has the following direct 
sum decomposition:
\begin{align}
\cA=\oplus_{\lambda\in \Lambda} \cA_\lambda. \label{P-Decomp}
\end{align}

%Furthermore, from the proof of the existence of 
%the Jordan-Chevalley decomposition for 
%linear maps of finite dimensional vector spaces 
%(e.g., see \cite[Proposition $4.2$]{H}) and 
%that for LF linear maps in \cite[Proposition $1.3.8$]{E} 
%it is easy to see that for each $\lambda\in \lambda$  
%\begin{align}
%\cA_\lambda=\sum_{i=1}^\infty \Ker (\lambda\I-\phi)^i
%\end{align}

Furthermore, each $\cA_\lambda$ $(\lambda\in \Lambda)$ is $\phi$ 
(and also $\Ps$ and $\Pn$) invariant, 
whence $\cA_\lambda\cA_\mu\subseteq \cA_{\lambda\mu}$ for all 
$\lambda, \mu\in \Lambda$. In particular, $\cA_1$ is a
$K$-subalgebra of $\cA$. Therefore we have
\begin{align}
\im\delta=\oplus_{\lambda\in \Lambda} \delta(\cA_\lambda). \label{ImP-Decomp}
\end{align}

\begin{lemma} \label{P-preLma(-1)}
Assume that $K$ is algebraically closed. Then
\begin{enumerate}
\item[$1)$] $\cA^\delta\subseteq \cA_1$.
\item[$2)$] $\im \delta=\Pn(\cA_1) \oplus 
\bigoplus_{1\ne \lambda\in \Lambda} \cA_\lambda$. 
\item[$3)$] $e\cA,\, \cA e  \subseteq \im \delta$ for all 
 idempotents $e\in \cA^\delta\cap\im \delta$.  
\end{enumerate}
\end{lemma}

\pf $1)$ Note that by the uniqueness of the Jordan-Chevalley decomposition of 
$\delta=\delta_{\frak s}+\delta_{\frak n}$ it is easy to see that 
$\delta_{\frak s}=\I-\Ps$ and $\delta_{\frak n}=-\Pn$.
Then by \cite[Proposition $1.3.9$, $i)$]{E} we have 
\begin{align} \label{P-preLma(-1)-peq0}
\cA^\delta=\Ker (\I-\Ps)\cap \Ker \Pn =\cA_1\cap \Ker \Pn.
\end{align}
Hence $\cA^\delta\subseteq \cA_1$. 
 
$2)$ Note first that for all $\lambda\in \Lambda$, we have 
$\delta(\cA_\lambda)\subseteq \cA_\lambda$ 
and 
\begin{align}
\delta \mid_{\cA_\lambda}=(\I_{\cA_\lambda}-\Ps\mid_{\cA_\lambda})-\Pn\mid_{\cA_\lambda}
=(1-\lambda) \I_{\cA_\lambda}-\Pn\mid_{\cA_\lambda}.\label{P-preLma(-1)-peq1}
\end{align}

Since $\Pn$ is LN, by Lemma \ref{InvLma}  
$\delta\mid_{\cA_\lambda}$ $(1\ne \lambda\in \Lambda)$ is a $K$-linear 
automorphism of $\cA_\lambda$. Hence 
$\cA_\lambda\subseteq \im \delta$ for 
all $1\ne \lambda\in \Lambda$. 
Furthermore, since $\delta\mid_{\cA_1}=-\Pn\mid_{\cA_1}$, by Eq.\,(\ref{ImP-Decomp})  
the statement follows. 

%$3)$ Note that $\Pn$ is a LN $K$-derivation of $\cA$ 
%(as pointed out above)   
%and $e\in \Ker \Dn$ by Eq.\,(\ref{D-preLma(-1)-peq1}). 
%Applying lemma \ref{LN-MainThm-1} to $\Dn$ we have 
%$e\cA,\, \cA e,\ \cA e\cA \subseteq \im \Dn$. 
%Therefore it suffices to show $\im\Dn\subseteq \im D$. 
%
%Since $\cA_\lambda$ is $\Dn$ invariant for all $\lambda\in \Lambda$, 
%we have 
%$$
%\im \Dn=\Dn(\cA_0)\oplus \bigoplus_{0\ne \lambda\in\Lambda} \Dn (\cA_\lambda).
%$$ 
%Since $\Ds (\cA_\lambda)\subseteq \cA_\lambda$ for all $\lambda\in \Lambda$, 
%by statement $2)$ we have $\im\Dn\subseteq \im D$, as desired. 

$3)$ By statement $1)$ we have $e\in \cA_1$, and    
by Statement $2)$,  $e\in \Pn(\cA_1)$. 
Since $(-\Pn)\mid_{\cA_1}=\delta\mid_{\cA_1}$ by Eq.\,(\ref{P-preLma(-1)-peq1}), 
we see that $-\Pn\mid_{\cA_1}$ is a LN $K$-$\cE$-derivation of $\cA_1$. Applying 
Lemma \ref{LN-MainThm-1} to the $K$-algebra $\cA_1$ and the $K$-$\cE$-derivation 
$(-\Pn)$ of $\cA_1$ we have 
$e\cA_1, \, \cA_1 e \subseteq \Pn(\cA_1)=\delta(\cA_1) \subseteq \im \delta$.  
Note that for all $1\ne \lambda \in \Lambda$, we also have  
$e\cA_\lambda, \, \cA_\lambda e \subseteq \cA_\lambda$ 
(for $e\in \cA_1$). Then by statement $2)$ we see that 
statement $3)$ follows.
\epfv

\begin{rmk}
From the proof of Lemma \ref{P-preLma(-1)}, $3)$ it is easy to see that we also have 
$\cA_\lambda e \cA_\mu \subseteq \im \delta$ for all 
idempotents $e\in \cA^\delta\cap\im \delta$ and all 
$\lambda,\mu\in \Lambda$ with $\lambda\mu\ne 1$. 
In particular, if $\lambda^{-1}\not \in \Lambda$ for all 
$0, 1\ne \lambda \in\Lambda$, then we also have $(e)\subseteq\im\delta$, 
as the $K$-derivation case of Theorem \ref{MainThm-1}. 
In general, it is still unknown weather or not $\cA_\lambda e \cA_{\lambda^{-1}} \subseteq \im \delta$ for all $\lambda\in\Lambda$ with $\lambda^{-1}\in \Lambda$.  
\end{rmk}

The rest of the proof of Theorem \ref{MainThm-1}, $2)$ for $\delta$, i.e., 
without assuming that $K$ is algebraically closed, 
is similar as that of the proof of Theorem \ref{MainThm-1}, $1)$  
at the end of the previous section. So we skip it here. 

Note also that Theorem \ref{MainThm-1}, $2)$ for the LN 
$K$-$\cE$-derivations of $\cA$ has been established 
in Lemma \ref{LN-MainThm-1} in the previous section.

%
%\begin{lemma} 
%Let $u\in \cA$ such that $E(u)=1$. 
%Then $u$ is transcendental over $K$.
%\end{lemma} 
%
%\pf Assume otherwise. Let $p(t)$ be the minimal polynomial of 
%$u$ over $K$.
%
%Since $1=E(u)=u-\phi(u)$, we have $\phi(u)=u-1$. 
%Then $0=\phi(p(u))=p(u-1)$. Since $\deg p(t-1)=\deg p(t)$ 
%and $p(t-1)$ is also monic, we have $p(t)=p(t-1)$. 
%Since char.\,$K$=0, this is impossible since otherwise, 
%$p(t)$ would have infinitely many distinct roots 
%in the algebraic closure $\bar K$ of $K$. 
%Contradiction.
%\epfv
% 

\section{\bf Proof of Proposition \ref{Surj-one}}\label{S5}

In this section we give a proof for Proposition \ref{Surj-one}. 
Note first that if $1 \in \cA^\delta$, 
then Proposition \ref{Surj-one} immediately follows  
from Theorem \ref{MainThm-1}. In particular, this is the case 
when $\delta$ is a $K$-derivation of $\cA$ (or $\delta=\I-\phi$ 
with $\phi(1)=1$). So we need only to show the proposition 
for LF (locally finite) $K$-$\cE$-derivations of $\cA$.

Throughout this section we let $\delta$ be a LF $K$-$\cE$-derivations 
of $\cA$ and write $\delta=\I-\phi$ with $e\!:=\phi(1)$. 
Note that $\phi$ is also LF and $e$ is an idempotent of $\cA$, 
for $\phi$ is a $K$-algebra endomorphism of $\cA$.  

\begin{lemma} \label{Lma-5.1}
For each $i\ge 0$, set $e_i=\phi^i(1)$. Then  
\begin{enumerate}
  \item[$1)$] $e_ie_j=e_je_i=e_j$ for all $0\le i\le j$. 
  \item[$2)$] there exists $d\ge 1$ such that $e_k=e_d$ for all $k\ge d$.   
\end{enumerate}
\end{lemma} 
 
\pf  $1)$ For all $0\le i\le j$, consider
$$
e_i e_j=\phi^i(1)\phi^{i}(\phi^{j-i}(1))=\phi^i\big(1\cdot\phi^{j-i}(1)\big)=
\phi^j(1)=e_j.
$$

Similarly, we also have $e_je_i=e_j$.
 
$2)$ If $e=0$, then $\phi(1)=0$, whence $\phi=0$ for 
$\phi\in\cEnd_K(\cA)$. In this case we may choose $d=1$. Assume $e\ne 0$ and 
let $V$ be the $K$-subspace spanned over $K$ by $\phi^i(1)$ $(i\ge 0)$. Then $V$ is 
$\phi$-invariant and finite dimensional over $K$, for $\phi$ is LF. 
Let $p(t)=t^d+\sum_{i=0}^{d-1}c_it^i\in K[t]$ be the minimal polynomial of $\phi\mid_V$.  
Then by applying $p(\phi)$ to $1$ we get

\begin{align}\label{Lma-5.1-peq1}
e_d+\sum_{i=0}^{d-1}c_i e_i=0.
\end{align}

Multiplying $e_d$ (from the left or the right) to the equation above 
and by statement $1)$ we get $(1+\sum_{i=0}c_i)e_d=0$.
If $e_d=0$, then the statement obviously holds. 
So we assume that $e_d\ne 0$ and get  
  
\begin{align}\label{Lma-5.1-peq2}
1+\sum_{i=0}^{d-1}c_i=0.
\end{align}
 
Let $0\le j\le d-1$ such that 
$c_j\ne 0$ but $c_i=0$ for all 
$j<i\le d-1$. Multiplying $e_j$  to Eq.\,(\ref{Lma-5.1-peq1})  
and by statement $1)$ and Eq.\,(\ref{Lma-5.1-peq2}) we get 
$$
0=e_d+(\sum_{i=0}^j c_i)e_j=e_d+(\sum_{i=0}^{d-1} c_i)e_j
=e_d-e_j.
$$ 

Hence $e_d=e_j$. By statement $1)$ again  
$e_d=e_d e_{d-1}=e_je_{d-1}=e_{d-1}$. Then by using the induction on $k\ge d$,  
it is easy to see that the statement indeed holds. 
\epfv 

\underline{\bf Proof of Proposition \ref{Surj-one}}:
As pointed out at the beginning of this section we need only 
to show the $K$-$\cE$-derivation case of the proposition.

Let $\delta$ be a LF $K$-$\cE$-derivation of $\cA$ with 
$1\in \im\delta$. Let $e=\phi(1)$,  
$d\ge 1$ and $e_d$ be as in Lemma \ref{Lma-5.1}. 
If $e_d=0$, i.e., $\phi^d(1)=0$, then  
$\phi^d=0$, for $\phi^d\in\cEnd_K(\cA)$. 
In other words, $\phi$ is (locally) nilpotent. 
Then by Lemma \ref{InvLma}, $\delta$ is invertible, 
i.e., the proposition holds in this case. 
 
Assume $e_d\ne 0$. Then by Lemma \ref{Lma-5.1}, 
$2)$ we have $\phi(e_d)=e_d$, whence $\phi^d(e_d)=e_d$. 
Let $\delta_d\!:= \I-\phi^d $. Then 
$\delta_d$ is also a LF $K$-$\cE$-derivation 
of $\cA$ (for $\phi^d\in\cEnd_K(\cA)$ and is LF) and 
$e_d\in \cA^{\delta_d}$. 

Next, we show that $e_d$ also lies in $\im \delta_d$. 
Let $u\in \cA$ such that $\delta(u)=1$. 
Then $\phi(u)=u-1$, and inductively we have 
\begin{align*}
\phi^d(u)=u-1-e_1-\cdots-e_{d-1}.
\end{align*}
Multiplying $e_d$ from the left to the equation above 
and by Lemma \ref{Lma-5.1}, $1)$ we have 
$$
\phi^d(e_d u)=e_d\phi^d(u)=e_du-de_d.
$$
Then $\phi^d(\frac{e_du}{d})=\frac{e_du}d-e_d$, which means 
$e_d=\delta_d(\frac{e_du}{d})\in \im\delta_d$. 
Therefore $e_d\in\cA^{\delta_d}\cap\, \im\delta_d$. 
Applying Theorem \ref{MainThm-1}, $2)$ to the $K$-$\cE$-derivation 
$\delta_d$ and the idempotent $e_d$ 
we have $e_d\cA\subseteq \im\delta_d$. 
 
One the other hand, since $\phi^d(1-e_d)=e_d-e_d=0$, we have 
$\phi^d ((1-e_d)\cA)=0$, i.e., the restriction of $\delta_d$ on 
$(1-e_d)\cA$ is the identity map of $(1-e_d)\cA$, whence 
$(1-e_d)\cA\subseteq \im\delta_d$. Since $\cA=e_d\cA+(1-e_d)\cA$, 
we have $\im\delta_d=\cA$. 
Furthermore, since $\delta_d=\I-\phi^d=(\I-\phi)\sum_{i=0}^{d-1}\phi^i$, 
we have $\im \delta_d \subseteq \im (\I-\phi)$, whence 
$\im\delta=\im (\I-\phi)=\cA$, as desired. 
\epfv
 
%
%
%\subsection{\bf Idempotent in the Image of an Ideal}
%
%\begin{lemma} 
%Assume that $D$ is a locally nilpotent derivation of $\cA$ with a slice $s$.
%Then 
%
%$1)$ $[s, \cA^D]\subseteq \cA^D$.
%
%$2)$ Assume that $\cA^D$ is reduced, i.e., does not have any nonzero nilpotent, 
%which is equivalent to that $a^2\ne 0$ for all nonzero 
%$a\in\cA^D$, then each idempotent $e$ of $\cA$ lies in $\cA^D$. 
%\end{lemma} 
%\pf $1)$ For all $a\in \cA^D$, consider
%$$
%D([s,a])=D(sa-as)=a-a=0.
%$$ 
%Hence $[s, a]\subset \cA^D$.
%
%$2)$ By lemma \ref{} we may write $e$ uniquely as 
%$e=\sum_{i=0}^d c_is^i$ for some $d\ge 0$ and 
%$c_i$'s in $\cA^D$ with $c_d s^d \ne 0$.
%
%If $d=0$, we are done. If $d\ge 1$, 
%then by plugging the expression into the equation 
%$e^2=e$ and using $1)$ and freeness of $s$ over $\cA^D$ 
%we have $c_d^2=0$, whence $c_d=0$. 
%Contradiction.
%\epfv
%


\begin{thebibliography}{FLM2}


\bibitem[BCW]{BCW} H. Bass, E. Connell and D. Wright, {\it The Jacobian Conjecture, Reduction of Degree and Formal Expansion of the Inverse}. Bull.  Amer. Math.  Soc.  \textbf{7}, (1982), 287--330. [MR 83k:14028]. 


\bibitem[B]{Bergen} J. Bergen, {\it Derivations in Prime Rings}. Canad. Math. Bull. {\bf 26} (1983), 267-270.

\bibitem[BFF]{BFF} M. Bre\v{s}ar, A. Fo\v{s}ner and M. Fo\v{s}ner, 
{\it A Kleinecke-Shirokov Type Condition with Jordan Automorphisms}.  
Studia Math. {\bf 147} (2001), no.\,3, 237-242. 

\bibitem[BV]{BV} M. Bre\v{s}ar and AR Villena, {\it The Noncommutative Singer–Wermer Conjecture and $\phi$-Derivations}. J. London Math. Soc. {\bf 66} (2002), 710-720.

\bibitem[DEZ]{DEZ} H. Derksen, A. van den Essen and W. Zhao, {\it The Gaussian Moments Conjecture and the Jacobian Conjecture}. To appear in {\it Israel J. Math.}. See also arXiv:1506.05192 [math.AC]. 

\bibitem[E1]{E0} A. van den Essen, {\it The Exponential Conjecture and the Nilpotency Subgroup of the Automorphism Group of a Polynomial Ring}. Prepublications. Univ. Aut\`onoma de Barcelona, April 1998. 

\bibitem[E2]{E} A. van den Essen, {\it Polynomial Automorphisms and the Jacobian Conjecture}.
Prog. Math., Vol.190, Birkh\"auser Verlag, Basel, 2000. %[MR1790619].

\bibitem[E3]{E2} A. van den Essen, {\it Introduction to Mathieu Subspaces}. 
``International Short-School/Conference on Affine Algebraic Geometry and the Jacobian Conjecture" at Chern Institute of Mathematics, Nankai University, Tianjin, China. 
July 14-25, 2014.
 
\bibitem[EH]{EH} A. van den Essen and L. C. van Hove, {\it Mathieu--Zhao Spaces}. 
             To appear. 

\bibitem[EN]{EN} A. van den Essen and S. Nieman, {\it Mathieu--Zhao Spaces of Univariate Polynomial Rings with Non-zero Strong Radical}.  J. Pure Appl. Algebra, 
{\bf 220} (2016), no.\,9, 3300--3306.


%\bibitem[E2]{E2} A. van den Essen,, {\it The Amazing Image Conjecture.} Preprint, arXiv:1006.5801v1[math.AG].

%\bibitem[K]{Kurosch} A. G. Kurosch, {\it Problems in Ring Theory Which Are Related to 
%the Burnside Problem for Periodic Groups}. Izv. Akad. Nauk SSSR 5, no. {\bf 3} (1941), 
%233–240.
%
%\bibitem[Br]{Bresar} M. Bre\v sar, {\it Semiderivations of Prime Rings}. 
%Proc. Amer. Math. Soc. 108 (1990), no. {\bf 4}, 859–860. [MR1007488] (90g:16030)

%Classified semi-der of prime rings.

%\bibitem[Ch]{Chuang}  Chen-Lian Chuang, 
%{\it On the structure of semiderivations in prime rings}.   
%Proc. Amer. Math. Soc. {\bf 108} (1990), 867-869. 
%
%%Classified semi-der of prime rings.

%\bibitem[Ch1]{Chang1} Jui Chi Chang, {\it On semiderivations of Prime Rings}. 
%Chinese J. Math. 12 (1984), no. {\bf 4}, 255–262. [MR0774289] (86g:16050)

% The author replaces derivations by semiderivations (introduced by J. Bergen [Canad. Math. Bull. 26 (1983), no. 3, 267–270; MR 85a:16079]), and proves several commutativity theorems, many of which are extensions of earlier results due to I. N. Herstein [ibid. 21 (1978), no. 3, 369–370; MR0506447 (58 #22182); ibid. 22 (1979), no. 4, 509–511; MR0563766 (81b:16025)]. 
%   Suppose f,f1,f2 are nonzero semiderivations of R. Any of the following conditions on the prime ring R will force R to be commutative: (i) f(R)\subset Z; (ii) [f(R),f(R)]=0 and chR\ne 2; (iii) [f(R),f(R)]\subset Z; (iv) f1f2(R)\subset Z; (v) there exists a nonzero element r in R such that rf(R)\subset Z; (vi) [r,f(r)]\subset Z for all r in R, and chR\ne 2.
 


\bibitem[EWZ]{EWZ} A. van den Essen, D. Wright and W. Zhao, 
{\it Images of Locally Finite Derivations of Polynomial Algebras in Two Variables}. J. Pure Appl. Algebra {\bf 215} (2011), no.9, 2130-2134. [MR2786603]. See also arXiv:1004.0521[math.AC].


\bibitem[GN]{GN} P. Gabriel and Y. Nouaz\'e, {\it Id\'eax Premiers de l'Alg\`ebra  
Enveloppante d'une Alg\`ebra de Lie Nilpotente}. J. Algebra {\bf 6} (1967), 77-99.

%\bibitem[H]{H} J. E. Humphreys, (1972), 
%{\it Introduction to Lie Algebras and Representation Theory}.  
%Graduate Texts in Mathematics, Springer, 1972. %ISBN 978-0-387-90053-7

\bibitem[J]{J} N. Jacobson, {\it Structure of Rings}. Amer. Math. Soc. Coll. Pub. {\bf 37}, 
Amer. Math. Soc. Providence R. I., 1956.

\bibitem[K]{K} O. H. Keller, 
{\it Ganze Gremona-Transformationen}. Monats. Math. Physik {\bf 47} (1939), no.\,1, 299-306. [MR1550818].


%\bibitem[Kh1]{Kharchenko1} V.K. Kharchenko, {\it Generalized identities with automorphisms}. Algebra Logika 14 (1975) 132–148. [MR0399153] (53 \#3004)
%
%\bibitem[Kh2]{Kharchenko2} V.K. Kharchenko, {\it Differential identities of prime rings}.  Algebra Logika 17 (1978) 220–238 (English transl, in Algebra Logic 17 (1978) 154–168). 
%MR0541758 (81f:16025)
%
%\bibitem[Kh3]{Kharchenko3} V.K. Kharchenko, {\it Automorphisms and Derivations of Associative Rings}. Math. Appl. (Soviet Ser.), vol. 69, Kluwer Academic Publishers Group, Dordrecht, 1991, translated from the Russian by L. Yuzina. MR1174740 (93i:16048)

 


\bibitem[M]{Ma} O. Mathieu, {\it Some Conjectures about Invariant Theory and Their Applications.} Alg\`ebre non commutative, groupes quantiques et invariants (Reims, 1995), 263--279, S\'emin. Congr., 2, Soc. Math. France, Paris, 1997. [MR1601155].
 


\bibitem[W]{W} D. Wright, {\it On the Jacobian Conjecture}. Illinois J. Math., 
{\bf 25} (1981), no.\,3, 423--440. 
 


%\bibitem[Z1]{HNP} W. Zhao, {\it Hessian Nilpotent Polynomials and the Jacobian Conjecture}, Trans. Amer. Math. Soc. {\bf 359} (2007), no. 1, 249--274 (electronic). [MR2247890]. See also math.CV/0409534.
%
%\bibitem[Z2]{GVC} W. Zhao, {\it A Vanishing Conjecture on Differential Operators with Constant Coefficients}, Acta Mathematica Vietnamica, vol 32 (2007), no. 3, 259--286.  [MR2368014]. See also arXiv:0704.1691v2 [math.CV].
%

\bibitem[Z1]{IC} W. Zhao, {\it Images of Commuting  Differential Operators of Order One with Constant Leading Coefficients}.  
J. Alg. {\bf 324} (2010),  no. 2, 231--247. [MR2651354].  
See also arXiv:0902.0210 [math.CV]. 

\bibitem[Z2]{GIC} W. Zhao, {\it Generalizations of the Image Conjecture and the Mathieu Conjecture}.  J. Pure Appl. Alg. {\bf 214} (2010), 1200-1216. See also arXiv:0902.0212 [math.CV].  
%
%\bibitem[Z2]{GMS} W. Zhao, {\it A Generalization of Mathieu Subspaces to Modules of Associative Algebras}. Preprint. 

\bibitem[Z3]{MS} W. Zhao, {\it Mathieu Subspaces of Associative Algebras}. J. Alg. 
{\bf 350} (2012), no.2, 245-272. [MR2859886]. See also arXiv:1005.4260 [math.RA].

\bibitem[Z4]{Open-LFNED} W. Zhao, {\it Some Open Problems on Locally Finite or Locally Nilpotent Derivations and $\cE$-Derivations}. Preprint. %See also arXiv:

%\bibitem[Z5]{Idem} W. Zhao, {\it Idempotents in Intersection of the Kernel and the Image of Locally Finite Derivations and $\cE$-derivations}. Preprint. %See also arXiv:

\bibitem[Z5]{Alg-LFNED} W. Zhao, {\it The LFED and LNED Conjectures 
for Algebraic Algebras}. 
Preprint. %See also arXiv:

\bibitem[Z6]{LaurentPolyCase} W. Zhao, {\it The LFED and LNED Conjectures 
for Laurent Polynomial Algebras}. 
Preprint. %See also arXiv:

\bibitem[Z7]{OneVariableCase} W. Zhao, {\it Images of Ideals under Derivations and $\cE$-Derivations of  Univariate Polynomial Algebras  over a Field  of Characteristic Zero}. Preprint. %See also arXiv:

\end{thebibliography}
\end{document}